\documentclass[10pt,a4paper]{amsart}

\usepackage{amssymb}
\usepackage{amsthm}
\usepackage{amsfonts}
\usepackage{microtype}
\usepackage{mathrsfs}
\usepackage{graphicx}
\usepackage{color}
\usepackage{latexsym}
\usepackage{rotating}
\usepackage{xspace}
\usepackage[all]{xy}
\usepackage{longtable}
\usepackage{leftidx}
\usepackage{mathtools}
\usepackage{titlesec}
\usepackage{tikz}
\usetikzlibrary{calc}
\usetikzlibrary{arrows}
\usetikzlibrary{decorations.markings}
\usepackage{geometry}

\newtheorem{theorem}{Theorem}[section]
\numberwithin{equation}{theorem}
\newtheorem{lemma}[theorem]{Lemma}

\newtheorem{corollary}[theorem]{Corollary}

\theoremstyle{definition}
\newtheorem{definition}[theorem]{Definition}
\newtheorem{example}[theorem]{Example}
\newtheorem{notation}[theorem]{Notation}
\newtheorem{remark}[theorem]{Remark}

\theoremstyle{conjecture}

\newtheorem{question}[theorem]{Question}

\newcommand{\Ass}{\operatorname{Ass}}
\newcommand{\im}{\operatorname{im}}

\newcommand{\grade}{\operatorname{grade}}

\newcommand{\Assh}{\operatorname{Assh}}
\newcommand{\Spec}{\operatorname{Spec}}

\newcommand{\ara}{\operatorname{ara}}

\newcommand{\cd}{\operatorname{cd}}

\newcommand{\Ht}{\operatorname{ht}}
\newcommand{\id}{\operatorname{id}}
\newcommand{\fd}{\operatorname{fd}}

\newcommand{\pd}{\operatorname{pd}}

\newcommand{\V}{\operatorname{V}}
\newcommand{\RH}{\operatorname{H}}
\newcommand{\RE}{\operatorname{E}}

\newcommand{\Ext}{\operatorname{Ext}}
\newcommand{\Supp}{\operatorname{Supp}}

\newcommand{\Tor}{\operatorname{Tor}}
\newcommand{\Hom}{\operatorname{Hom}}

\newcommand{\Ann}{\operatorname{Ann}}
\newcommand{\Rad}{\operatorname{Rad}}

\newcommand{\Vdim}{\operatorname{Vdim}}
\newcommand{\depth}{\operatorname{depth}}

\newcommand{\Max}{\operatorname{Max}}

\newcommand{\vpl}{\operatornamewithlimits{\varprojlim}}

\newcommand{\lo}{\longrightarrow}
\newcommand{\fm}{\frak{m}}
\newcommand{\fp}{\frak{p}}
\newcommand{\fq}{\frak{q}}
\newcommand{\fa}{\frak{a}}
\newcommand{\fb}{\frak{b}}

\newcommand{\fn}{\frak{n}}

\newcommand{\suchthat}{\;\ifnum\currentgrouptype=16 \middle\fi|\;}

\newenvironment{prf}[1][Proof]{\begin{proof}[\bf #1]}{\end{proof}}

\makeatletter
\newcommand{\holim@}[2]{%
\vtop{\m@th\ialign{##\cr
\hfil$#1\operator@font holim$\hfil\cr
\noalign{\nointerlineskip\kern1.5\ex@}#2\cr
\noalign{\nointerlineskip\kern-\ex@}\cr}}%
}
\newcommand{\holim}{%
\mathop{\mathpalette\holim@{\rightarrowfill@\textstyle}}\nmlimits@
}
\makeatother

\makeatletter
\def\@secnumfont{\bfseries}
\def\section{\@startsection{section}{1}%
\z@{.7\linespacing\@plus\linespacing}{.5\linespacing}%
{\normalfont\Large\bfseries\filcenter}}
\def\subsection{\@startsection{subsection}{2}%
\z@{.5\linespacing\@plus.7\linespacing}{-.5em}%
{\normalfont\large\bfseries}}
\makeatother
\begin{document}

\author[P. Pourghobadian, K. Divaani-Aazar and A. Rahimi]
{Parisa Pourghobadian, Kamran Divaani-Aazar and Ahad Rahimi}

\title[Dualities and equivalences of ...]
{Dualities and equivalences of the category of relative Cohen-Macaulay modules}

\address{P. Pourghobadian, Department of Mathematics, Faculty of Mathematical Sciences, Alzahra
University, Tehran, Iran.}
\email{paparpourghobadian@gmail.com}

\address{K. Divaani-Aazar, Department of Mathematics, Faculty of Mathematical Sciences, Alzahra
University, Tehran, Iran.}
\email{kdivaani@ipm.ir}

\address{A. Rahimi, Department of Mathematics, Razi University, Kermanshah, Iran.}
\email{ahad.rahimi@razi.ac.ir}



\begin{abstract} In this paper, we establish the global analogues of some dualities and equivalences in local algebra by developing
the theory of relative Cohen-Macaulay modules. Let $R$ be a commutative Noetherian ring (not necessarily local) with identity and
$\fa$ a proper ideal of $R$. The notions of $\fa$-relative dualizing modules and $\fa$-relative big Cohen-Macaulay modules are
introduced. With the help of $\fa$-relative dualizing modules, we establish the global analogue of the duality on the subcategory
of Cohen-Macaulay modules in local algebra. Lastly, we investigate the behavior of the subcategory of $\fa$-relative Cohen-Macaulay
modules and $\fa$-relative generalized Cohen-Macaulay modules under Foxby equivalence.
\end{abstract}

\maketitle

\tableofcontents

\section{Introduction}

Throughout this article, the word ring refers to commutative Noetherian rings with identity. Let $\fa$ be a proper ideal of a ring
$R$. Our objective is to prove the global analogues of some dualities and equivalences in local algebra. To accomplish this, we
extend the theory of relative Cohen-Macaulay modules by introducing and studying the concepts of relative dualizing modules and
relative big Cohen-Macaulay modules. Due to their application in this paper, these concepts deserve further exploration.

Hellus and Schenzel \cite{HSc} came up with the notion of cohomologically complete intersection ideals in 2008. By their definition,
for a Gorenstein local ring $(R,\fm)$, the ideal $\fa$ is {\it cohomologically complete intersection} if $\RH_{\fa}^{i}\left(R\right)
=0$ for all $i\neq \Ht \fa$. Rahro Zargar and Zakeri \cite{RZ2} define a finitely generated $R$-module $M$ to be $\fa$-{\it relative
Cohen-Macaulay} if $\RH_{\fa}^{i}\left(M\right)=0$ for all $i\neq \cd(\fa,M)$, where $\cd(\fa,M)$ is the largest integer $i$
for which $\RH_{\fa}^{i}\left(M\right)\neq 0$. Obviously, when $R$ is a Gorenstein local ring, $\fa$ is cohomologically complete
intersection if and only if $R$ is $\fa$-relative Cohen-Macaulay. Several authors followed up on relative Cohen-Macaulay modules;
see e.g. \cite{HeSt1, Sc2, R, JR, Ra2, RZ1, CH, Ra1, DGTZ1, Sc1, DGTZ2}. In particular, in  \cite{DGTZ2} the authors introduce the
notion of $\fa$-relative system of parameters which appears to be quite useful in studying relative Cohen-Macaulay modules. The
paper proceeds as follows.

Section 2 is devoted to the basic properties of relative system of parameters and hyperhomology.

We introduce the notion of relative dualizing modules in Section 3. They can be characterized by the vanishing of certain Ext
modules; see Lemma \ref{3.5} and Remark \ref{3.5a}. Despite the fact that relative dualizing modules are not necessarily finitely
generated, our results reveal that they inherit many properties of dualizing modules; see Theorems \ref{3.4}, \ref {3.7},
\ref{3.9} and \ref{3.91}.  If $T$ is a module-finite $R$-algebra such that $R$ is $\fa$-relative Cohen-Macaulay and $T$ is
$\fa T$-relative Cohen-Macaulay, then Theorem \ref{3.7} asserts that $\Omega_{\fa T}\cong \Ext_R^{t}(T,\Omega_{\fa})$, where $t:=\cd\left(\fa,R\right)-\cd\left(\fa T,T\right)$, and for a $\fb$-relative Cohen-Macaulay ring $S$, $\Omega_{\fb}$ denotes
its $\fb$-relative dualizing module. The main ingredient of the proof of Theorem \ref{3.7} is Lemma \ref{3.6}, which expresses
the minimal injective cogenerator of a module-finite $R$-algebra according to the minimal injective cogenerator of $R$.

The theory of big Cohen-Macaulay modules and algebras is a central topic in commutative algebra. Hochster's conjecture on the
existence of big Cohen-Macaulay algebras remained unsolved for over four decades. In Section 4, we define the notion of relative
big Cohen-Macaulay modules. We show that relative big Cohen-Macaulay modules behave very much like the ordinary big Cohen-Macaulay
modules; see Theorems \ref{3.10} and \ref{3.12}.

Section 5 deals with dualities. There is a well-known duality of categories:
\begin{displaymath}
\xymatrix{\text{CM}^n(R) \ar@<0.7ex>[rrr]^-{\Ext_R^{\dim R-n}(-,\omega_R)} &
{} & {} &  \text{CM}^n(R),  \ar@<0.7ex>[lll]^-{\Ext_R^{\dim R-n}(-,\omega_R)}}
\end{displaymath}
where $(R,\fm)$ is a Cohen-Macaulay local ring with a dualizing module $\omega_R$ and $\text{CM}^n(R)$ denotes the full subcategory
of $n$-dimensional Cohen-Macaulay $R$-modules. In \cite{CH}, Celikbas and Holm establish an analogue of the above duality for local
$\fa$-relative Cohen-Macaulay rings. Assuming that $R$ is an $\fa$-relative Cohen-Macaulay ring (not necessarily local), we prove
a relative analogue of the above duality; see Theorem \ref{4.1}. As a result, if $\mathfrak{J}$ denotes the Jacobson radical of
a complete semi-local ring $R$ and $R$ is $\mathfrak{J}$-relative Cohen-Macaulay, then for every integer $0\leq n\leq \dim R$, there
is a duality of categories:
\begin{displaymath}
\xymatrix{\text{CM}^n_{\mathfrak{J}}(R) \ar@<0.7ex>[rrr]^-{\Ext_R^{\dim R-n}(-,\Omega_{\mathfrak{J}})} &
{} & {} &  \text{CM}^n_{\mathfrak{J}}(R).  \ar@<0.7ex>[lll]^-{\Ext_R^{\dim R-n}(-,\Omega_{\mathfrak{J}})}}
\end{displaymath}
Here, for a proper ideal $\fb$ of a ring $S$, $\text{CM}^n_{\fb}(S)$ denotes the full subcategory of $\fb$-relative Cohen-Macaulay
$S$-modules with $\cd(\fb,M)=n$. For $\mathfrak{J}$-relative generalized Cohen-Macaulay modules, it is natural to ask whether the
above duality holds. In this regard, we show that if $(R,\fm)$ is a Cohen-Macaulay local ring with a dualizing module $\omega_R$
and $M$ a generalized Cohen-Macaulay $R$-module,  then the $R$-module $\Ext_R^{\dim R-\dim_RM}(M,\omega_{R})$ is also generalized
Cohen-Macaulay; see Theorem \ref{4.3}. We conclude Section 5 by proving a relative analogue of Grothendieck's local duality theorem;
see Theorem \ref{4.4}.

Let $C$ be a semidualizing module of $R$, and $\mathscr{A}_C\left(R\right)$ and $\mathscr{B}_C\left(R\right)$ denote, respectively,
the Auslander and the Bass classes of $R$ with respect to $C$. Then, by Foxby equivalence,  there is an equivalence of categories:
\begin{displaymath}
\xymatrix{\mathscr{A}_C\left(R\right) \ar@<0.7ex>[rrr]^-{C\otimes_R-} &
{} & {} & \mathscr{B}_C\left(R\right).  \ar@<0.7ex>[lll]^-{\Hom_R\left(C,-\right)}}
\end{displaymath}
Let $n$ be a non-negative integer, and $\text{gCM}^n_{\fa}(R)$ denote the full subcategory of $\fa$-relative generalized Cohen-Macaulay $R$-modules $M$ with $\cd(\fa,M)=n$. In Section 6, we prove the following equivalence of categories:
\begin{displaymath}
\xymatrix{\mathscr{A}_C\left(R\right)\bigcap \text{CM}^n_{\fa}(R) \ar@<0.7ex>[rrr]^-{C\otimes_R-} &
{} & {} & \mathscr{B}_C\left(R\right)\bigcap \text{CM}^n_{\fa}(R),  \ar@<0.7ex>[lll]^-{\Hom_R\left(C,-\right)}}
\end{displaymath}
and
\begin{displaymath}
\xymatrix{\mathscr{A}_C\left(R\right)\bigcap \text{gCM}^n_{\fa}(R) \ar@<0.7ex>[rrr]^-{C\otimes_R-} &
{} & {} & \mathscr{B}_C\left(R\right)\bigcap \text{gCM}^n_{\fa}(R).  \ar@<0.7ex>[lll]^-{\Hom_R\left(C,-\right)}}
\end{displaymath}

\section{Prerequisites}

In this section, we provide some background material on relative system of parameters and hyperhomology which will be used in the
rest of this work.

\vspace{.3cm}
${\bf Relative\ system\ of\ parameters:}$
\vspace{.5cm}

Below, we recall some definitions and results from \cite{DGTZ2}. This paper discusses, among other things, local cohomology modules $$\RH_{\fa}^{i}\left(M\right):=\varinjlim \limits_{n\in \mathbb{N}} \text{Ext}_R^i\left(R/\fa^n,M\right); \  i\in \mathbb{N}_0.$$
We denote by $\cd(\fa,M)$ the {\it cohomological dimension of $M$ with respect to $\fa$} which is the largest integer $i$ such that $\RH_{\fa}^{i}\left(M\right)\neq 0$.

\begin{definition}\label{2.1} Let $M$ be a finitely generated $R$-module and $\fa$ an ideal of $R$ with $M\neq \fa M$.
\begin{enumerate}
\item[(i)] Let $c:=\cd\left(\fa,M\right)$. A sequence $x_{1}, x_2, \ldots, x_{c}\in \fa$ is called $\fa$-{\it relative system of
parameters}, $\fa$-s.o.p, of $M$ if $$\Rad\left(\langle x_{1}, x_2, \ldots, x_{c}\rangle+\Ann_{R}M\right)=\Rad\left(\fa+\Ann_{R}M
\right).$$
\item[(ii)] {\it Arithmetic rank} of $\fa$ with respect to $M$, $\ara\left(\fa,M \right)$, is defined as the infimum of the integers
$n\in \mathbb{N}_0$ such that there exist $x_1, x_2, \ldots, x_n\in R$ satisfying $$\Rad\left(\langle x_{1}, x_2, \ldots, x_{n}\rangle +\Ann_{R}M\right)=\Rad\left(\fa+\Ann_{R}M\right).$$
\end{enumerate}
\end{definition}

Clearly, if $x_{1},x_2,\ldots, x_{c}\in R$ is an $\fa$-s.o.p of $M$, then for all $t_{1},\ldots,t_{c}\in \mathbb{N}$, every permutation
of $x_{1}^{t_{1}},\ldots, x_{c}^{t_{c}}$ is also an $\fa$-s.o.p of $M$. One may easily check that $\cd\left(\fa,M\right)\leq \ara
\left(\fa,M \right)$. Obviously, $\ara\left(\fa,R\right)=\ara\left(\fa\right)$.

Over a local ring, every finitely generated module possesses a system of parameters, but this is not true for $\fa$-relative systems
of parameters. In this regard, we have:

\begin{lemma}\label{2.2} (See \cite[Lemma 2.2]{DGTZ2}.) Let $M$ be a finitely generated $R$-module and $\fa$ an ideal of $R$ with
$M\neq \fa M$. Then $\fa$ contains an $\fa$-s.o.p of $M$ if and only if $\cd\left(\fa,M\right)=\ara\left(\fa,M \right)$.
\end{lemma}

\begin{theorem}\label{2.3} (See \cite[Lemma 2.4 and Theorem 2.7]{DGTZ2}.) Let $\fa$ be an ideal of $R$, $M$ a finitely generated
$R$-module with $M\neq \fa M$ and $c=\cd\left(\fa,M\right)$. Assume that $\cd\left(\fa,M\right)=\ara\left(\fa,M \right)$ and
$x_1,\ldots, x_{c}\in \fa$. Consider the following conditions:
\begin{enumerate}
\item[(i)] $x_1,\ldots, x_{c}$ is an $\fa$-s.o.p of $M$.
\item[(ii)] $\cd\left(\fa,M/\langle x_{1},x_2, \ldots, x_{i}\rangle M\right)=c-i$ for every $i=1, 2,\ldots, c$.
\end{enumerate}
Then (i) implies (ii). Furthermore, if $\fa$ is contained in the Jacobson radical of $R$, then (i) and (ii) are equivalent.
\end{theorem}

\begin{definition}\label{2.4} Let $\fa$ a proper ideal of $R$ and $M$ a finitely generated $R$-module. Then $M$ is said to be
$\fa$-{\it relative Cohen-Macaulay} if either $M=\fa M$ or $M\neq \fa M$ and $\grade(\fa,M)=\cd(\fa,M)$.
\end{definition}

\begin{theorem}\label{2.5} (See \cite[Theorem 3.3]{DGTZ2}.) Let $M$ be a finitely generated $R$-module and $\fa$ an ideal of $R$ with $\cd\left(\fa,M\right)=\ara\left(\fa,M \right)$. Consider the following conditions:
\begin{enumerate}
\item[(i)] $M$ is $\fa$-relative Cohen-Macaulay.
\item[(ii)] Every $\fa$-s.o.p of $M$ is an $M$-regular sequence.
\item[(iii)] There exists an $\fa$-s.o.p of $M$ which is an $M$-regular sequence.
\end{enumerate}
Then (i) and (iii) are equivalent. Furthermore, if $\fa$ is contained in the Jacobson radical of $R$, then all three conditions
are equivalent.
\end{theorem}

\vspace{.3cm}
${\bf Derived\ Hom\ and\ tensor\ bifunctors:}$
\vspace{.5cm}

Throughout, the symbol $\simeq$ will refer to isomorphisms in the category $\mathcal{D}(R)$, the derived category of $R$-modules. An
object in $\mathcal{D}(R)$ is an $R$-complex $X$ displayed in the standard homological style $$X= \cdots \rightarrow X_{i+1}
\xrightarrow {\partial^{X}_{i+1}} X_{i} \xrightarrow {\partial^{X}_{i}} X_{i-1} \rightarrow \cdots.$$ Refer to \cite{AF}, \cite{C}
and \cite{Li} for more information.

We let $\mathcal{D}_{\sqsubset}(R)$ (res. $\mathcal{D}_{\sqsupset}(R)$) denote the full subcategory of $\mathcal{D}(R)$ consisting
of $R$-complexes $X$ with $\RH_{i}(X)=0$ for $i \gg 0$ (res. $i \ll 0$), and $\mathcal{D}_{\square}(R):=\mathcal{D}_{\sqsubset}(R)\cap \mathcal{D}_{\sqsupset}(R)$. Moreover, let $\mathcal{D}^{f}(R)$ denote the full subcategory of $\mathcal{D}(R)$ consisting of
$R$-complexes $X$ with finitely generated homology modules. In addition, we feel free to use any combination of the subscripts and
the superscript as in $\mathcal{D}^{f}_{\square}(R)$, with the obvious meaning of the intersection of the two subcategories involved.

An $R$-complex $P$ of projective modules is said to be semi-projective if the functor $\Hom_{R}(P,-)$ preserves quasi-isomorphisms.
By a semi-projective resolution of an $R$-complex $X$, we mean a quasi-isomorphism $P\xrightarrow {\simeq} X$ in which $P$ is a
semi-projective $R$-complex. Dually, an $R$-complex $I$ of injective modules is said to be semi-injective if the functor
$\Hom_{R}(-,I)$ preserves quasi-isomorphisms. By a semi-injective resolution of an $R$-complex $X$, we mean a quasi-isomorphism
$X\xrightarrow {\simeq} I$ in which $I$ is a semi-injective $R$-complex. Semi-projective and semi-injective resolutions exist for
any $R$-complex. Moreover, any right-bounded $R$-complex of projective modules is semi-projective, and any left-bounded $R$-complex
of injective modules is semi-injective.

Let $X$ and $Y$ be two $R$-complexes. Then ${\bf R}\Hom_{R}(X,Y)$ can be computed by $${\bf R}\Hom_{R}(X,Y)\simeq \Hom_{R}(P,Y) \simeq \Hom_{R}(X,I),$$ where $P\xrightarrow {\simeq} X$ is any semi-projective resolution of $X$, and $Y\xrightarrow {\simeq} I$ is any
semi-injective resolution of $Y.$  Moreover, we let $$\Ext^{i}_{R}(X,Y):=\RH_{-i}({\bf R}\Hom_{R}(X,Y))$$ for every $i\in
\mathbb{Z}$. Likewise, $X\otimes_{R}^{\bf L}Y$ can be computed by $$X\otimes_{R}^{\bf L}Y \simeq P\otimes_{R}Y \simeq X\otimes_{R}Q,$$
where $P\xrightarrow {\simeq} X$ is any semi-projective resolution of $X$, and $Q\xrightarrow {\simeq} Y$ is any semi-projective
resolution of $Y.$ Moreover, we let $$\Tor^{R}_{i}(X,Y):=\RH_{i}(X\otimes_{R}^{\bf L}Y)$$ for every $i \in \mathbb{Z}$.

For any complexes $X,Y$ and $Z$, there are the following natural isomorphisms in $\mathcal{D}(R)$.
\begin{enumerate}
\item[]  {\bf Shifts:} Let $i,j$ be two integers. Then
\begin{equation}
\Sigma^{i}X\otimes_R^{{\bf L}}\Sigma^{j}Y \simeq \Sigma^{j+i}(X\otimes_R^{{\bf
L}}Y)
\end{equation} and
\begin{equation}
{\bf R}\Hom_R(\Sigma^{i}X,\Sigma^{j}Y)\simeq \Sigma^{j-i}{\bf R}\Hom_R(X,Y).
\end{equation}
\item[] {\bf Adjointness:} If $X\in \mathcal{D}_{\sqsupset}(R)$, $Y\in \mathcal{D}(R)$ and $Z\in \mathcal{D}_{\sqsubset}(R)$, then
\begin{equation}
{\bf R}\Hom_R(X\otimes_R^{{\bf L}}Y,Z)\simeq {\bf R}\Hom_R(X,{\bf R}\Hom_R(Y,Z)).
\end{equation}
\item[] {\bf Tensor evaluation:} Assume that $X\in \mathcal{D}^f_{\sqsupset}(R)$, $Y\in \mathcal{D}_{\square}(R)$ and
$Z\in \mathcal{D}_{\sqsupset}(R)$. If either the projective dimension of $X$ or the flat dimension of $Z$ is finite, then
\begin{equation}
{\bf R}\Hom_R(X,Y)\otimes_R^{{\bf L}}Z\simeq {\bf R}\Hom_R(X,Y\otimes_R^{{\bf L}}Z).
\end{equation}
\item[] {\bf Hom evaluation:} Assume that $X\in \mathcal{D}^f_{\sqsupset}(R)$, $Y\in \mathcal{D}_{\square}(R)$ and
$Z\in \mathcal{D}_{\sqsubset}(R)$. If either the projective dimension of $X$ or the injective dimension of $Z$ is
finite, then
\begin{equation}
X\otimes_R^{{\bf L}}{\bf R}\Hom_R(Y,Z)\simeq {\bf R}\Hom_R({\bf R}\Hom_R(X,Y),Z).
\end{equation}
\end{enumerate}

In what follows, we will use the above isomorphisms without further explanation.

Let $X\in \mathcal{D}(R)$ and $s$ be an integer. If $\RH_i(X)=0$ for all $i\neq s$, then $X\simeq \Sigma^{s}\RH_s(X)$.

\vspace{.5cm}
${\bf Derived\ local\ cohomology:}$
\vspace{.5cm}

Let $\fa$ be an ideal of $R$. The right derived functor of the $\fa$-section functor $$\Gamma_{\fa}(-):=\underset{n}{\varinjlim}
\ \Hom_R(R/\fa^n,-)$$ exists in $\mathcal{D}(R)$. Let $X\in\mathcal{D}(R)$. Then the complex ${\bf R}\Gamma_{\fa}(X)$ is defined by
${\bf R}\Gamma_{\fa}(X) :=\Gamma_{\fa}(I)$, where $I$ is any semi-injective resolution of $X$. For any integer $i$, the $i$-th
local cohomology module of $X$ with respect to $\fa$ is defined by $\RH^i_{\fa}(X):=\RH_{-i}({\bf R}\Gamma_{\fa}(X))$. Let
$\Check{C}(\underline{\fa})$ denote the \v{C}ech complex of $R$ on a set $\underline{\fa}=a_1,\ldots, a_n$ of generators of $\fa$.
By \cite[Theorem 1.1(iv)]{Sc3},
\begin{equation}
{\bf R}\Gamma_{\fa}(X)\simeq X\otimes_R^{{\bf L}} \Check{C}(\underline{\fa})\simeq X
\otimes_R^{{\bf L}} {\bf R}\Gamma_{\fa}(R).
\end{equation}

For every $X\in \mathcal{D}(R)$ and any integer $s$, by $(2.5.6)$, we may conclude that
\begin{equation}
{\bf R}\Gamma_{\fa}(\Sigma^{s}X)\simeq
\Sigma^{s}{\bf R}\Gamma_{\fa}(X).
\end{equation}
As $\Check{C}(\underline{\fa})$ is a bounded complex of flat $R$-modules, tensor evaluation
along with $(2.5.6)$ yield the following known result:

\begin{lemma}\label{2.6} Let $\fa$ be an ideal of $R$, $X\in \mathcal{D}^f_{\sqsupset}(R)$ and $Y\in \mathcal{D}_{\square}(R)$.
Then $${\bf R}\Gamma_{\fa}({\bf R}\Hom_R(X,Y))\simeq {\bf R}\Hom_R(X,{\bf R}\Gamma_{\fa}(Y)).$$
\end{lemma}

\section{Relative dualizing modules}

We begin this section by defining relative dualizing modules. In what follows, $\Max R$ denotes for the set of maximal ideals of $R$.

\begin{definition}\label{3.1} Let $\fa$ be a proper ideal of $R$ and $c=\cd(\fa,R)$. Assume that $R$ is $\fa$-relative
Cohen-Macaulay. The $\fa$-{\it relative dualizing module} of $R$ is defined by $$\Omega_{\fa}:=\Hom_R(\RH_{\fa}^c(R),
\underset{\fm\in \Max R}\bigoplus \RE_R(R/\fm)).$$
\end{definition}

The dualizing module of a complete Cohen-Macaulay local ring $(R,\fm)$ is the same as its $\fm$-relative dualizing module.
A local ring $(R,\fm)$ admits a dualizing module if and only if it is Cohen-Macaulay and a quotient of a Gorenstein local
ring. Nevertheless,  by the definition, $R$ admits an $\fa$-relative dualizing module if and only if it is $\fa$-relative
Cohen-Macaulay.

The following two lemmas are required to prove Theorem \ref{3.4}.

\begin{lemma}\label{3.3} (See \cite[Proposition 2.1]{RZ2}.) Let $\fa$ be an ideal of $R$ and $n$ a non-negative integer. Let $M$ be
an $R$-module such that $\RH^{i}_{\fa}\left(M\right)=0$ for all $i\neq n$. Then $$\Ext_R^i(R/\fa,\RH_{\fa}^n(M))\cong
\Ext_R^{i+n}(R/\fa,M)$$ for all $i\geq 0$.
\end{lemma}

\begin{lemma}\label{3.5} Let $\fa$ be a proper ideal of $R$, $\RE:=\underset{\fm\in \Max R}\bigoplus \RE_R(R/\fm)$ and $0\leq n\leq c=\cd(\fa,R)$
an integer. Let $M$ be an $R$-module such that $\RH_{\fa}^i(M)=0$ for all $i\neq n$. Assume that $R$ is $\fa$-relative Cohen-Macaulay. Then
\begin{enumerate}
\item[(i)] $\Ext_R^i(M,\Omega_{\fa})=
\begin{cases} \Hom_R(\RH_{\fa}^n(M),\RE) \hspace{.5cm} \text{if} \hspace{.5cm}    i=c-n \\
0 \hspace{3.1cm}   \text{if} \hspace{.5cm}   i\neq c-n.
\end{cases}
$
\item[(ii)] $\Ext_R^i(R/\fa,\Omega_{\fa})=
\begin{cases} \Hom_R(R/\fa,\RE) \hspace{.5cm}   \text{if} \hspace{.5cm}   i=c \\
0 \hspace{2.6cm}   \text{if} \hspace{.5cm}    i\neq c.
\end{cases}
$
\end{enumerate}
\end{lemma}

\begin{prf} (i) Observe that
$$\begin{array}{ll}
{\bf R}\Hom_R\left(M,\Omega_{\fa}\right)&\simeq {\bf R}\Hom_R\left(M,{\bf R}\Hom_R(\RH_{\fa}^c(R),\RE)\right)\\
&\simeq {\bf R}\Hom_R\left(M,{\bf R}\Hom_R(\Sigma^{c}{\bf R}\Gamma_{\fa}(R),\RE)\right)\\
&\simeq \Sigma^{-c}{\bf R}\Hom_R\left(M,{\bf R}\Hom_R({\bf R}\Gamma_{\fa}(R),\RE)\right)\\
&\simeq \Sigma^{-c}{\bf R}\Hom_R\left(M\otimes_R^{{\bf L}}{\bf R}\Gamma_{\fa}(R),\RE\right)\\
&\simeq \Sigma^{-c}{\bf R}\Hom_R\left({\bf R}\Gamma_{\fa}(M),\RE\right)\\
&\simeq \Sigma^{-c}{\bf R}\Hom_R\left(\Sigma^{-n}\RH_{\fa}^n(M),\RE\right)\\
&\simeq \Sigma^{-(c-n)}{\bf R}\Hom_R\left(\RH_{\fa}^n(M),\RE\right)\\
&\simeq \Sigma^{-(c-n)}\Hom_R\left(\RH_{\fa}^n(M),\RE\right).\\
\end{array}$$
The first and last steps in the above display hold by definition. As $\RH_i({\bf R}\Gamma_{\fa}(M))=0$ for all $i\neq -n$, it 
follows that ${\bf R}\Gamma_{\fa}(M)\overset{(\dagger)}\simeq \Sigma^{-n}\RH_{\fa}^n(M)$. Similarly, since $R$ is $\fa$-relative Cohen-Macaulay, we have ${\bf R}\Gamma_{\fa}(R)\simeq \Sigma^{-c}\RH_{\fa}^c(R)$, which implies $\RH_{\fa}^c(R)\overset{(\ddagger)}\simeq  \Sigma^{c}{\bf R}\Gamma_{\fa}(R)$. Clearly, $(\dagger)$ and $(\ddagger)$ give rise to the sixth and second quasi-isomorphisms, 
respectively. The third and seventh steps follow by $(2.5.2)$. The fourth step is by $(2.5.3)$. Finally, the fifth step is by 
$(2.5.6)$. Consequently, $${\bf R}\Hom_R\left(M,\Omega_{\fa}\right)\simeq \Sigma^{-(c-n)}\Hom_R\left(\RH_{\fa}^n(M),\RE\right),$$ 
which by taking homology yields (i).

(ii) As $R/\fa$ is $\fa$-torsion, it turns out that $\RH_{\fa}^0(R/\fa)=R/\fa$ and $\RH_{\fa}^i(R/\fa)=0$ for all $i\neq 0$.
Thus, $R/\fa$ is $\fa$-relative Cohen-Macaulay and $n:=\cd(\fa,R/\fa)=0$, and so (i) yields the claim.
\end{prf}

\begin{remark}\label{3.5a} Let $R$ be an $\fa$-relative Cohen-Macaulay ring and $c=\cd(\fa,R)$. Lemma \ref{3.5}(i) provides a
characterization of the $\fa$-relative dualizing module. Let $\Delta$ be an $R$-module such that the assertion of Lemma \ref{3.5}(i)
holds if we replace $\Delta$ by $\Omega_{\fa}$. Then $$\Delta\cong \Hom_R(R,\Delta)\cong \Hom_R\left(\RH_{\fa}^c(R),\RE\right)
=\Omega_{\fa}.$$
\end{remark}

\begin{theorem}\label{3.4} Let $\fa$ be a proper ideal of $R$, $\RE=\underset{\fm\in \Max R}\bigoplus \RE_R(R/\fm)$ and $c=\cd(\fa,R)$.
Assume that $R$ is $\fa$-relative Cohen-Macaulay. Then
\begin{enumerate}
\item[(i)] $\id_R(\Omega_{\fa})=c$.
\item[(ii)] $\RH_{\fa}^i(\Omega_{\fa})=
\begin{cases} \Gamma_{\fa}(\RE)  \hspace{.4cm}  \text{if} \hspace{.4cm}  i=c \\
0 \hspace{1.1cm} \text{if} \hspace{.4cm} i\neq c.
\end{cases}
$
\item[(iii)] ${\bf R}\Hom_R\left(\Omega_{\fa},\Omega_{\fa}\right)\simeq \underset{\fm\in \V(\fa)\cap \Max R}\prod \widehat{R_{\fm}}$.
\end{enumerate}
\end{theorem}

\begin{prf} (i) Set $(-)^\vee:=\Hom_R(-,\RE)$. So, $\Omega_{\fa}=\RH_{\fa}^c(R)^{\vee}$. By \cite[Theorem 3.2]{RZ1}, it is known that $\fd_R(\RH_{\fa}^c(R))=c$. As $(-)^\vee$ is a faithfully exact contravariant functor on the category of $R$-modules, it turns out that $\id_R(\Omega_{\fa})=c$.

(ii) As $\id_R(\Omega_{\fa})=c$, it follows that $\RH_{\fa}^i(\Omega_{\fa})=0$ for all $i>c$. Lemma \ref{3.5}(ii) implies that $\Ext_R^{i}(R/\fa,\Omega_{\fa})=0$ for all $i<c$. So, \cite[Proposition 5.3.15]{St} yields that $\RH_{\fa}^i(\Omega_{\fa})=0$ for all $i<c$. Consequently, $\RH_{\fa}^i(\Omega_{\fa})=0$ for all $i\neq c$. As $\id_R(\Omega_{\fa})=c$, the $R$-module $\Omega_{\fa}$ admits an injective resolution
$$0\lo I_0\lo \cdots \lo I_{-c+1}\xrightarrow{d_{-c+1}} I_{-c}\lo 0.$$ Then, $\RH_{\fa}^c(\Omega_{\fa})=\Gamma_{\fa}( I_{-c})/\im(\Gamma_{\fa}(d_{-c+1}))$. It is known that $\Gamma_{\fa}(I)$ is injective for every injective $R$-module $I$. Therefore, as $\RH_{\fa}^i(\Omega_{\fa})=0$ for all $i\neq c$, the sequence $$0\lo \Gamma_{\fa}(I_0)\lo \cdots \lo \Gamma_{\fa}(I_{-c+1})\lo \im(\Gamma_{\fa}(d_{-c+1}))\lo 0$$ is split exact. In particular, the $R$-module $\im(\Gamma_{\fa}(d_{-c+1}))$ is injective. Thus, the exact sequence $$0\lo \im(\Gamma_{\fa}(d_{-c+1}))\lo \Gamma_{\fa}( I_{-c})\lo \RH_{\fa}^c(\Omega_{\fa})\lo 0$$  splits, and so the $R$-module $\RH_{\fa}^c(\Omega_{\fa})$ is injective.

Next, Lemma \ref{3.3} and Lemma \ref{3.5}(ii) imply that
\begin{equation}
\Hom_R(R/\fa,\RH_{\fa}^c(\Omega_{\fa}))\cong \Ext_R^c(R/\fa,\Omega_{\fa})\cong (R/\fa)^\vee.
\end{equation}
For every injective $R$-module $I$, by Matlis' theorem, one has $$I\cong
\underset{\fp \in \Ass_RI}\bigoplus \RE_R(R/\fp)^{\mu^0(\fp,I)},$$ where $\mu^0(\fp,I)$ is the $0$-th Bass number of $I$ with respect to
the prime ideal $\fp$.

For simplicity, we set $Y:=\RH_{\fa}^c(\Omega_{\fa})$ and $\overline{R}:=R/\fa$. As $Y$ is $\fa$-torsion, it follows that
$\Ass_RY\subseteq \V(\fa)$. For every prime ideal $\fp$ of $R$, one has $$\Hom_R(R/\fa,\RE_R(R/\fp))=
\begin{cases} \RE_{\overline{R}}(R/\fp)  \hspace{.4cm}  \text{if} \hspace{.4cm}  \fa\subseteq \fp\\
0 \hspace{1.6cm}  \text{if} \hspace{.4cm} \fa \nsubseteq \fp.
\end{cases}
$$
Hence, $$\Hom_R(R/\fa,Y)=\underset{\fp \in \Ass_RY}\bigoplus \RE_{\overline{R}}(R/\fp)^{\mu^0(\fp,Y)}$$ and $$(R/\fa)^\vee=
\underset{\fm\in \Max R\cap \V(\fa)}\bigoplus \RE_{\overline{R}}(R/\fm).$$ From this, by $(3.5.1)$, we conclude that $\Ass_RY=\Max R
\cap \V(\fa)$, and $\mu^0(\fm,Y)=1$ for all $\fm\in \Ass_RY$. Therefore, $$Y\cong \underset{\fm\in \Max R\cap \V(\fa)}
\bigoplus \RE_R(R/\fm)=\Gamma_{\fa}(\RE).$$

(iii) Set $X:=\V(\fa)\cap \Max R$. For every $\fm\in \Max R$, we can see that $$\Hom_R\left(\RE_R(R/\fm),\RE_R(R/\fm)\right)\cong
\Hom_{R_{\fm}}\left(\RE_{R_{\fm}}(R_{\fm}/\fm R_{\fm}),\RE_{R_{\fm}}(R_{\fm}/\fm R_{\fm})\right)\cong \widehat{R_{\fm}},$$ and
$$\Hom_R(\RE_R(R/\fm),\underset{\fn\in \Max R\setminus \{\fm\}}\bigoplus \RE_R(R/\fn))=0.$$ Therefore, (ii) and Lemma \ref{3.5}(i)
imply that
$$\begin{array}{ll}
{\bf R}\Hom_R(\Omega_{\fa},\Omega_{\fa})&\simeq  \Hom_R(\Gamma_{\fa}(\RE),\RE)\\
&\simeq \Hom_R\left(\underset{\fm\in X}\bigoplus \RE_R(R/\fm),\underset{\fm\in \Max R}\bigoplus \RE_R(R/\fm)\right)\\
&\simeq \underset{\fm\in X}\prod\Hom_R\left(\RE_R(R/\fm),\underset{\fm\in \Max R}\bigoplus \RE_R(R/\fm)\right)\\
&\simeq \underset{\fm\in X}\prod\Hom_R(\RE_R(R/\fm),\RE_R(R/\fm))\\
&\simeq \underset{\fm\in X}\prod \widehat{R_{\fm}}.
\end{array}$$
\end{prf}

Let $X$ be as in the proof of Theorem \ref{3.4}(iii). Then the ring $\underset{\fm\in X}\prod \widehat{R_{\fm}}$ is Noetherian if and
only if the set $X$ is finite. If the set $X$ is finite, then the ring $\underset{\fm\in X}\prod \widehat{R_{\fm}}$ is the $J$-adic
completion of $R$, where $J$ is the intersection of all maximal ideals of $R$ containing $\fa$.

It is known that a Cohen-Macaulay local ring $(R,\fm)$ with a dualizing module $\omega_R$ is Gorenstein if and only if $R\cong \omega_R$.
By Theorem \ref{3.4}(i), if $R$ possesses a proper ideal $\fa$ such that $R$ is $\fa$-relative Cohen-Macaulay and $R\cong \Omega_{\fa}$,
then $R$ is Gorenstein.

The next result expresses the minimal injective cogenerator of a module-finite $R$-algebra according to the minimal injective
cogenerator of $R$.

\begin{lemma}\label{3.6} Let $\phi: R\lo T$ be a module-finite ring homomorphism. Then there is a natural $T$-isomorphism
$$\Hom_R(T,\underset{\fm\in \Max R}\bigoplus \RE_R(R/\fm))\cong \underset{\fn\in \Max T}\bigoplus \RE_T(T/\fn).$$
\end{lemma}

\begin{prf} Set $X:=\{\phi^{-1}(\fn) \mid \fn\in \Max T\}$. Let $\fm\in \Max R$ and $W_{\fm}:=\Hom_R(T,\RE_R(R/\fm))$.
We claim that $W_{\fm}\neq 0$ if and only if $\fm\in X$.

First, assume that $\fm\in X$. Then, $\fm=\phi^{-1}(\fn)$ for some $\fn\in \Max T$. The inclusion $\Ann_TT=0\subseteq \fn$
yields that $\Ann_RT=\phi^{-1}(\Ann_TT)\subseteq \fm$. Hence, $$\fm\in \Supp_RT\cap \{\fm\}=\Supp_RT\cap \Ass_R(\RE_R(R/\fm))
=\Ass_RW_{\fm},$$ and so $W_{\fm}\neq 0$.

Next, suppose that $W_{\fm}\neq 0$. Then $\fm\in \Ass_RW_{\fm}$, and so $\fm=\phi^{-1}(\fn)$ for some $\fn\in \Ass_TW_{\fm}$.
As $\phi$ is an integral ring homomorphism, it follows that $\fn\in \Max T$, and so $\fm\in X$.

Now, let $\fm\in X$, and set $Y_{\fm}:=\{\fn\in \Max T \mid \phi^{-1}(\fn)=\fm\}$. Let $\fn\in Y_{\fm}$, and set $l_{\fn}:=T/\fn$.
In what follows, we denote the dimension of a vector space $V$ over a field $F$ by $\Vdim_FV$. As $l_{\fn}$ is a finite-dimensional
extension of the field $R/\fm$, we see that $e_{\fn}:=\Vdim_{R/\fm}l_{\fn}$ is finite. One has
$$
\begin{array}{ll}
\Hom_T(l_{\fn},W_{\fm})&\cong \Hom_T\left(l_{\fn},\Hom_R(T,\RE_R(R/\fm))\right)\\
&\cong \Hom_R\left(l_{\fn}\otimes_TT,\RE_R(R/\fm)\right)\\
&\cong \Hom_R\left(l_{\fn},\RE_R(R/\fm)\right)\\
&\cong \Hom_R\left(\overset{e_{\fn}}\bigoplus R/\fm,\RE_R(R/\fm)\right)\\
&\cong \overset{e_{\fn}}\bigoplus \Hom_R\left(R/\fm,\RE_R(R/\fm)\right)\\
&\cong \overset{e_{\fn}}\bigoplus R/\fm\\
&\cong l_{\fn},
\end{array}
$$
and so $$\mu^0(\fn,W_{\fm})=\Vdim_{l_{\fn}}\left(\Hom_T\left(l_{\fn},W_{\fm}\right)\right)=1.$$ In particular, it turns out that $\Ass_TW_{\fm}=Y_{\fm}$. As $W_{\fm}$ is an Artinian injective $T$-module,  by Matlis' theorem, we conclude that $W_{\fm}\cong
\underset{\fn\in Y_{\fm}} \bigoplus  \RE_T(T/\fn)$. Now, we have
$$
\begin{array}{ll}
\Hom_R(T,\underset{\fm\in \Max R}\bigoplus \RE_R(R/\fm))&\cong \underset{\fm\in \Max R}\bigoplus W_{\fm}\\
&\cong \underset{\fm\in X}\bigoplus W_{\fm}\\
&\cong \underset{\fm\in X}\bigoplus (\underset{\fn\in Y_{\fm}}\bigoplus  \RE_T(T/\fn))\\
&\cong \underset{\fn\in \Max T}\bigoplus \RE_T(T/\fn).
\end{array}
$$
\end{prf}

Let $T$ be a module-finite $R$-algebra and $\fa$ a proper ideal of $R$ such that $R$ is $\fa$-relative Cohen-Macaulay and
$T$ is $\fa T$-relative Cohen-Macaulay.  The next result asserts that $\Omega_{\fa T}\cong \Ext_R^{t}(T,\Omega_{\fa})$, where
$t:=\cd\left(\fa,R\right)-\cd\left(\fa T,T\right).$

\begin{theorem}\label{3.7} Let $\phi: R\lo T$ be a module-finite ring homomorphism. Let $\fa$ be a proper ideal of
$R$, and set $c=\cd\left(\fa,R\right)$ and $c':=\cd\left(\fa T,T\right)$. Assume that $R$ is $\fa$-relative Cohen-Macaulay
and $T$ is $\fa T$-relative Cohen-Macaulay. Then $\Omega_{\fa T}\cong \Ext_R^{c-c'}(T,\Omega_{\fa})$.
\end{theorem}

\begin{prf} Set $\RE=\underset{\fm\in \Max R}\bigoplus \RE_R(R/\fm)$. By Lemma \ref{3.5}(i) and Lemma \ref{3.6}, we have
$$
\begin{array}{ll}
\Ext_R^{c-c'}(T,\Omega_{\fa})&\cong \Hom_R(\RH_{\fa}^{c'}(T),\RE)\\
&\cong \Hom_R\left(\RH_{\fa T}^{c'}(T),\RE\right)\\
&\cong \Hom_R\left(\RH_{\fa T}^{c'}(T)\otimes_TT,\RE\right)\\
&\cong \Hom_T\left(\RH_{\fa T}^{c'}(T),\Hom_R(T,\RE)\right)\\
&\cong \Hom_T(\RH_{\fa T}^{c'}(T),\underset{\fn\in \Max T}\bigoplus \RE_T(T/\fn))\\
&\cong \Omega_{\fa T},
\end{array}
$$
as required.
\end{prf}

In many ways, relative dualizing modules behave like dualizing modules. This is illustrated by the next result.

\begin{theorem}\label{3.9} Let $\fa$ be a proper ideal of $R$. Assume that $R$ is $\fa$-relative
Cohen-Macaulay.
\begin{enumerate}
\item[(i)] $\Omega_{\fa (R/\langle \underline{x} \rangle)}\cong \Omega_{\fa}/\langle \underline{x} \rangle \Omega_{\fa}$ for every
$R$-regular sequence $\underline{x}=x_1, \dots, x_r\in \fa$, which is a part of an $\fa$-s.o.p of $R$.
\item[(ii)] If $(R,\fm,k)$ is a local ring, then $\Omega_{\fa \widehat{R}}\cong \Omega_{\fa}$.
\end{enumerate}
\end{theorem}

\begin{prf} (i) Let $\underline{x}=x_1, \dots, x_r\in \fa$ be an $R$-regular sequence, which is a part of an $\fa$-s.o.p of $R$.
By induction on $r$, we may and do assume that $r=1$. Set $c=\cd(\fa,R)$, $\overline{R}:=R/\langle x_1 \rangle$ and
$\RE=\underset{\fm\in \Max R}
\bigoplus \RE_R(R/\fm)$. Then, by Theorem \ref{2.3}, one has $$\cd(\fa \overline{R},\overline{R})=\cd(\fa,\overline{R})=c-1.$$
Consider the short exact sequence
\begin{equation}
0\lo R\overset{x_1}\lo R\lo \overline{R}\lo 0.
\end{equation}
It yields the exact sequence $$\cdots \lo \RH_{\fa}^i(R)\lo \RH_{\fa}^i(\overline{R})\lo  \RH_{\fa}^{i+1}(R)\lo \cdots,$$ from
which we obtain that $\RH_{\fa}^i(\overline{R})=0$ for all $i<c-1$. Thus, the ring $\overline{R}$ is $\fa \overline{R}$-relative
Cohen-Macaulay. The short exact sequence $(3.8.1)$ also implies the following exact sequence $$0\lo \Hom_R(\overline{R},\Omega_{\fa})\lo \Omega_{\fa}
\overset{x_1}\lo \Omega_{\fa}\lo \Ext_R^1(\overline{R},\Omega_{\fa})\lo 0,$$ and so $\Ext_R^1(\overline{R},\Omega_{\fa})\cong
\Omega_{\fa}/\langle x_1 \rangle \Omega_{\fa}$. Now, Theorem \ref{3.7} completes the argument.

(ii) First note that by the Flat Base Change theorem \cite[Theorem 4.3.2]{BS}, one has
$$\RH_{\fa \widehat{R}}^i(\widehat{R})\cong \RH_{\fa}^i(R)\otimes_R\widehat{R}$$ for all $i\geq 0$. This along with the faithful flatness
of $\widehat{R}$ yield that $\widehat{R}$ is $\fa \widehat{R}$-relative Cohen-Macaulay and $\cd(\fa \widehat{R},\widehat{R})=\cd(\fa,R)$.
Now, the claim follows from the following display:
$$\begin{array}{ll}
\Omega_{\fa \widehat{R}}&=\Hom_{\widehat{R}}\left(\RH_{\fa \widehat{R}}^c(\widehat{R}),\RE_{\widehat{R}}(k)\right)\\
&\cong \Hom_{\widehat{R}}\left(\RH_{\fa}^c(R)\otimes_R\widehat{R},\RE_{\widehat{R}}(k)\right)\\
&\cong \Hom_R\left(\RH_{\fa}^c(R),\Hom_{\widehat{R}}(\widehat{R},\RE_{\widehat{R}}(k))\right)\\
&\cong \Hom_R\left(\RH_{\fa}^c(R),\RE_{\widehat{R}}(k)\right)\\
&\cong \Hom_R\left(\RH_{\fa}^c(R),\RE_R(k)\right)\\
&\cong \Omega_{\fa} .
\end{array}$$
\end{prf}

\section{Relative big Cohen-Macaulay modules}

We start by defining relative big Cohen-Macaulay modules.

\begin{definition}\label{3.8} Let $\fa$ be a proper ideal of $R$ with $\cd(\fa,R)=\ara(\fa)$. Let $M$ be a (not necessarily finitely
generated) $R$-module.
\begin{enumerate}
\item[(i)] We say $M$ is $\fa$-{\it relative big Cohen-Macaulay}, if there exists an $\fa$-s.o.p of $R$ which is an $M$-regular
sequence.
\item[(ii)] We say $M$ is $\fa$-{\it relative balanced big Cohen-Macaulay}, if every $\fa$-s.o.p of $R$ is an $M$-regular sequence.
\end{enumerate}
\end{definition}

Let $\fa$ be a proper ideal of $R$ with $\cd(\fa,R)=\ara(\fa)$. If $R$ is $\fa$-relative Cohen-Macaulay, then by Theorem \ref{2.5}, it
is $\fa$-relative big Cohen-Macaulay. It is also evident that if $\cd(\fa,R)=0$, then every $R$-module is vacuously $\fa$-relatively
balanced big Cohen-Macaulay.

One may guess that the condition $\cd(\fa,R)=\ara(\fa)$ in the above definition holds for $\fa$-relative Cohen-Macaulay rings.
However, the following example shows that this is not the case.

\begin{example}\label{3.888} Let $\Bbbk$ be a field and $S=\Bbbk[[x,y,z,w]]$. Consider the elements $f=xw-yz$, $g=y^{3}-x^{2}z$,
and $h=z^{3}-y^2w$ of $S$. Let $R=S/\langle f \rangle$, and $\mathfrak{a}=\langle f,g,h \rangle/\langle f \rangle$. Then $R$ is
a local complete intersection ring of dimension $3$, $\cd(\mathfrak{a},R)=1$, and $\ara(\mathfrak{a})\geq 2$; see
\cite[Remark 2.1(ii)]{HeSt2}. On the other hand, since ${\RH}_{\fa}^0(R)=0$, it follows that $R$ is $\fa$-relative Cohen-Macaulay.
\end{example}

Next, we present a large class of $\fa$-relatively big Cohen-Macaulay modules, which aren't necessarily finitely generated.

\begin{theorem}\label{3.91} Let $\fa$ be a proper ideal of $R$ with $\cd(\fa,R)=\ara(\fa)$. Assume that $R$ is $\fa$-relative
Cohen-Macaulay. Then $\Omega_{\fa}$ is $\fa$-relative big Cohen-Macaulay, and, furthermore, if $\fa$ is contained in the Jacobson
radical of $R$, then $\Omega_{\fa}$ is $\fa$-relative balanced big Cohen-Macaulay.
\end{theorem}

\begin{prf} Set $c=\cd(\fa,R)$. The claim vacuously holds for $c=0$, so we may assume that $c\geq 1$. By Theorem \ref{2.5}, there is
an $\fa$-s.o.p $x_1, \dots, x_c$ of $R$ which is an $R$-regular sequence. By induction on $c$, we show that $x_1, \dots, x_c$ is an $\Omega_{\fa}$-regular sequence. As we saw in the proof of Theorem \ref{3.9}, the ring $\overline{R}:=R/\langle x_1 \rangle$ is $\fa \overline{R}$-relative Cohen-Macaulay and $\cd(\fa \overline{R},\overline{R})=c-1$. Set $\RE=\underset{\fm\in \Max R}\bigoplus
\RE_R(R/\fm)$. From the short exact sequence $$0\lo R\overset{x_1}\lo R\lo \overline{R}\lo
0,$$ we get the short exact sequence $$0\lo \RH_{\fa}^{c-1}(\overline{R})\lo \RH_{\fa}^{c}(R)\overset{x_1}\lo\RH_{\fa}^{c}(R)\lo 0.$$
Applying the contravariant exact functor $\Hom_R(-,\RE)$ on the latter exact sequence, yields the following exact sequence  $$0\lo\Omega_{\fa}\overset{x_1}\lo \Omega_{\fa}\lo \Hom_R(\RH_{\fa}^{c-1}(\overline{R}),\RE)\lo 0.$$
Hence, $x_1$ is a non-zero divisor on $\Omega_{\fa}$. This yields the claim for $c=1$. Next, suppose that $c>1$ and the claim is true for
$c-1$. As $x_2+\langle x_1 \rangle, \dots, x_c+\langle x_1 \rangle$ is an $\fa \overline{R}$-s.o.p of $\overline{R}$ and $x_2+\langle x_1
\rangle, \dots, x_c+\langle x_1 \rangle$ is an $\overline{R}$-regular sequence, by induction hypothesis, it follows that $x_2+\langle x_1
\rangle, \dots, x_c+\langle x_1 \rangle$ is an $\Omega_{\fa \overline{R}}$-regular sequence. But, by Theorem \ref{3.9}(i), we have
$\Omega_{\fa \overline{R}}\cong \Omega_{\fa}/\langle x_1 \rangle \Omega_{\fa}$, and so  $x_1, \dots, x_c$ is an $\Omega_{\fa}$-regular
sequence. Thus, $\Omega_{\fa}$ is $\fa$-relative big Cohen-Macaulay.

If $\fa$ is contained in the Jacobson radical of $R$, then by Theorem \ref{2.5} every $\fa$-s.o.p of $R$ is an $R$-regular
sequence, and so the above argument yields that $\Omega_{\fa}$ is $\fa$-relative balanced big Cohen-Macaulay.
\end{prf}

Next, we find that relative big Cohen-Macaulay modules inherit many properties of the ordinary big Cohen-Macaulay modules.

\begin{theorem}\label{3.10} Let $\fa$ be a proper ideal of $R$ with $\cd(\fa,R)=\ara(\fa)$, and set $c=\cd(\fa,R)$. Let $M$ be an
$\fa$-relative big Cohen-Macaulay $R$-module.
\begin{enumerate}
\item[(i)] $M/\langle x_1, \dots, x_n\rangle M$ is an $\fa(R/\langle x_1, \dots, x_n\rangle)$-relative big Cohen-Macaulay module
for every $\fa$-s.o.p $x_1, \dots, x_c$ of $R$ which is an $M$-regular sequence and all $1\leq n\leq c$.
\item[(ii)] If $M$ is balanced, then $M/\langle x_1, \dots, x_n\rangle M$ is an $\fa(R/\langle x_1, \dots, x_n\rangle)$-relative
balanced big Cohen-Macaulay module for all $\fa$-s.o.p $x_1, \dots, x_c$ of $R$ and all $1\leq n\leq c$.
\item[(iii)] $S^{-1}M$ is an $\fa S^{-1}R$-relative big Cohen-Macaulay module for every multiplicatively closed subset $S$ of
$R$ such that $S^{-1}(\RH_{\fa}^c(R))\neq 0$.
\item[(iv)] If $(R,\fm)$ is a local ring, then $M\otimes_R\widehat{R}$ is an $\fa\widehat{R}$-relative big Cohen-Macaulay module.
\end{enumerate}
\end{theorem}

\begin{prf} (i) Let $x_1, \dots, x_c$ be an $\fa$-s.o.p of $R$ which is an $M$-regular sequence and let $1\leq n\leq c$ be an integer.
Set $\underline{x}:=x_1, \dots, x_n$ and $\overline{R}:=R/\langle \underline{x} \rangle$. If $n=c$, then $\cd(\fa \overline{R},\overline{R})=\ara(\fa \overline{R})=0$, and so every $\overline{R}$-module is vacuously
$\fa \overline{R}$-relatively balanced big Cohen-Macaulay. So, we may and do assume that $n<c$. Then, $x_{n+1}+\langle \underline{x} \rangle,
\dots, x_c+\langle \underline{x} \rangle$ is an $\fa \overline{R}$-s.o.p of $\overline{R}$ and $x_{n+1}+\langle \underline{x} \rangle,
\dots, x_c+\langle \underline{x} \rangle$ is an $M/\langle \underline{x} \rangle M$-regular sequence. Hence,  $M/\langle \underline{x}
\rangle M$ is $\fa \overline{R}$-relative big Cohen-Macaulay.

(ii) Let $x_1, \dots, x_c$ be an $\fa$-s.o.p of $R$ and $1\leq n\leq c$ an integer. Set $\underline{x}:=x_1, \dots, x_n$ and $\overline{R}:=R/\langle \underline{x} \rangle$. Similar to (i), we can suppose that $n<c$.  Let $y_{n+1}+\langle \underline{x} \rangle,
\dots, y_c+\langle \underline{x}\rangle$ be an $\fa \overline{R}$-s.o.p of $\overline{R}$. It is easy to verify that $x_1, \dots, x_n, y_{n+1}, \dots, y_c$ is
an $\fa$-s.o.p of $R$. As $M$ is balanced, $x_1, \dots, x_n, y_{n+1}, \dots, y_c$ is an $M$-regular sequence. Hence, $y_{n+1}+
\langle \underline{x} \rangle, \dots, y_c+\langle \underline{x}\rangle$ is an $M/\langle \underline{x} \rangle M$-regular sequence,
and so $M/\langle \underline{x} \rangle M$ is $\fa \overline{R}$-relative balanced big Cohen-Macaulay.

(iii) Since $S^{-1}(\RH_{\fa}^i(R))\cong \RH_{\fa S^{-1}R}^i(S^{-1}R)$ for all $i\geq 0$, and $S^{-1}(\RH_{\fa}^c(R))
\neq 0$, it follows that $$c\leq \cd(\fa S^{-1}R,S^{-1}R)\leq \cd(\fa,R)=c.$$ Consequently, $\cd(\fa S^{-1}R,S^{-1}R)=c$. By our
assumption, there is an
$\fa$-s.o.p $x_1, \dots, x_c$ of $R$ which is an $M$-regular sequence. The equality $$\Rad(\langle x_1, \dots, x_c\rangle_R)=\Rad(\fa)$$
implies that $$\Rad(\langle x_1/1, \dots, x_c/1\rangle_{S^{-1}R})=\Rad(\fa S^{-1}R).$$ Thus, $ x_1/1, \dots, x_c/1$ is an
$\fa S^{-1}R$-s.o.p of $S^{-1}R$. Also, from the fact that $x_1, \dots, x_c$ is an $M$-regular sequence, we readily conclude that $x_1/1,
\dots, x_c/1$ is an $S^{-1}M$-regular sequence. Therefore, $S^{-1}M$ is $\fa S^{-1}R$-relative big Cohen-Macaulay.

(iv) As we saw in the proof of Theorem \ref{3.9}(ii), one has $\cd(\fa \widehat{R},\widehat{R})=\cd(\fa,R)$. There is an $\fa$-s.o.p
$x_1, \dots, x_c$ of $R$ which is an $M$-regular sequence. Let $\psi: R\lo \widehat{R}$ be the natural ring monomorphism. Similar to
the proof of (iii), we can deduce that $$\Rad(\langle \psi(x_1), \dots, \psi(x_c)\rangle_{\widehat{R}})=\Rad(\fa \widehat{R}),$$ and
so $\psi(x_1), \dots, \psi(x_c)$ is an $\fa\widehat{R}$-s.o.p of $\widehat{R}$. Moreover, it is easy to verify that $\psi(x_1), \dots,
\psi(x_c)$ is an $M\otimes_R\widehat{R}$-regular sequence.  Therefore, $M\otimes_R\widehat{R}$ is $\fa \widehat{R}$-relative big
Cohen-Macaulay.
\end{prf}

The following lemma is needed to prove the last theorem of this section.

\begin{lemma}\label{3.11} Let $\fa$ be a proper ideal of $R$ with $\cd(\fa,R)=\ara(\fa)$, and set $c=\cd(\fa,R)$. Let $M$ be an
$\fa$-relative big Cohen-Macaulay $R$-module. Then
$$\RH_{\fa}^i(M)=
\begin{cases} \RH_{\fa}^c(R)\otimes_RM \hspace{.5cm}   \text{if} \hspace{.5cm}  i=c \\
0 \hspace{2.35cm}   \text{if} \hspace{.5cm}  i\neq c.
\end{cases}
$$
\end{lemma}

\begin{prf} There is an $\fa$-s.o.p $x_1, \dots, x_c$ of $R$ which is an $M$-regular sequence. By \cite[Propositions 5.3.5 and 5.3.15]{St},
it turns out that $\RH_{\fa}^i(M)=0$ for all $i<c$. Since $\cd(\fa,M)\leq c$, it follows that $\RH_{\fa}^i(M)=0$ for all $i>c$.
Finally, as the covariant functor $\RH_{\fa}^c(-)$ is right exact and commutes with direct limits, it follows that $\RH_{\fa}^c(M)\cong \RH_{\fa}^c(R)\otimes_RM$.
\end{prf}

Next, over a relative Cohen-Macaulay ring, we construct a relative big Cohen-Macaulay module from a  given one.

\begin{theorem}\label{3.12} Let $\fa$ be an ideal of $R$ contained in its Jacobson radical. Assume that $R$ is $\fa$-relative
Cohen-Macaulay and $\cd(\fa,R)=\ara(\fa)$. Then, for every $\fa$-relative (balanced) big Cohen-Macaulay $R$-module $M$, the $R$-module $\Hom_R(M,\Omega_{\fa})$ is also $\fa$-relative (balanced) big Cohen-Macaulay.
\end{theorem}

\begin{prf} First, let $M$ be an $\fa$-relative big Cohen-Macaulay $R$-module. Set $c=\cd(\fa,R)$ and $\mathcal{M}:=\Hom_R(M,
\Omega_{\fa})$. Obviously, we may assume that $c\geq 1$. There is an $\fa$-s.o.p $x_1,\dots, x_c$ of $R$ which is an $M$-regular
sequence. As $\fa$ is contained in the Jacobson radical of $R$, Theorem \ref{3.91} implies that $x_1,\dots, x_c$ is an
$\Omega_{\fa}$-regular sequence. By induction on $c$, we show that $x_1, \dots, x_c$ is an $\mathcal{M}$-regular sequence. Let
$c=1$. Applying the functor $\Hom_R(M,-)$ on the short exact sequence $$0\lo \Omega_{\fa}\overset{x_1}\lo \Omega_{\fa}\lo
\Omega_{\fa}/\langle x_1 \rangle \Omega_{\fa}\lo 0,$$ yields the exact sequence
\begin{equation}
0\lo \mathcal{M}\overset{x_1}\lo \mathcal{M}\lo \Hom_R(M,\Omega_{\fa}/\langle x_1\rangle \Omega_{\fa})\lo 0.
\end{equation}
Note that $\Ext_R^1(M,\Omega_{\fa})=0$ by Lemma \ref{3.11}
and Lemma \ref{3.5}(i). In particular, $x_1$ is a non-zero divisor on $\mathcal{M}$, and so the claim holds in this case.

Next, suppose that $c>1$ and the claim holds for $c-1$. Set $\overline{R}:=R/\langle x_1 \rangle$ and $\overline{M}:=M/\langle x_1
\rangle M$. Then $\fa \overline{R}$ is contained in the Jacobson radical of $\overline{R}$ and $\cd(\fa \overline{R},\overline{R})
=\ara(\fa \overline{R})$. On the other hand, $\overline{R}$ is $\fa \overline{R}$-relative Cohen-Macaulay by
\cite[Corollary 3.5]{DGTZ2}. Clearly, $x_2+\langle x_1 \rangle,\dots, x_c+\langle x_1 \rangle$ is both an $\fa \overline{R}$-s.o.p of $\overline{R}$ and an $\overline{M}$-regular
sequence. Hence, the $\overline{R}$-module $\overline{M}$ is $\fa \overline{R}$-relative big Cohen-Macaulay. By induction hypothesis,
it follows that $x_2+\langle x_1 \rangle,\dots, x_c+\langle x_1 \rangle$ is an $\Hom_{\overline{R}}(\overline{M},\Omega_{\fa \overline{R}})$-regular sequence.
Now,  by adjointness, Theorem \ref{3.9}(i) and $(4.6.1)$, we have:
$$\begin{array}{ll}
\Hom_{\overline{R}}(\overline{M},\Omega_{\fa \overline{R}})&\cong \Hom_{\overline{R}}\left(M\otimes_R\overline{R},\Omega_{\fa \overline{R}}\right)\\
&\cong \Hom_R\left(M,\Hom_{\overline{R}}(\overline{R},\Omega_{\fa \overline{R}})\right)\\
&\cong \Hom_R\left(M,\Omega_{\fa \overline{R}}\right)\\
&\cong \Hom_R\left(M,\Omega_{\fa}/\langle x_1 \rangle \Omega_{\fa}\right)\\
&\cong \mathcal{M}/\langle x_1 \rangle \mathcal{M}.
\end{array}
$$
Therefore,  $x_1,\dots, x_c$ is an $\mathcal{M}$-regular sequence.

Now, let $M$ be an $\fa$-relative balanced big Cohen-Macaulay $R$-module. Then, the above argument shows that the $R$-module $\Hom_R(M,\Omega_{\fa})$ is $\fa$-relative balanced big Cohen-Macaulay.
\end{prf}

We end this section by proposing the following question.

\begin{question}\label{3.13} Let $\fa$ be a proper ideal of $R$ such that $\ara(\fa)=\cd(\fa,R)$. Does it follow that $R$ admits an
$\fa$-relative big Cohen-Macaulay module?
\end{question}

\section{Dualities}

Let $\fa$ be a proper ideal of $R$. We call a finitely generated $R$-module $M$ $\fa$-{\it relative maximal Cohen-Macaulay} if
$\grade(\fa,M)=\cd(\fa,R)$. Also, a finitely generated $R$-module $M$ is called $\fa$-{\it relative generalized Cohen-Macaulay}
if $\RH_{\fa}^i(M)$ is finitely generated for all $i<\cd(\fa,M)$; see \cite[Definition 3.2]{DGTZ1}.

\begin{notation}\label{4.0}  Let $\fa$ be a proper ideal of $R$ and $n$ a non-negative integer. In what follows:
\begin{itemize}
\item[(i)] $\text{MCM}_{\fa}(R)$ stands for the full subcategory of $\fa$-relative maximal Cohen-Macaulay $R$-modules.
\item[(ii)] $\text{CM}^n_{\fa}(R)$ stands for the full subcategory of $\fa$-relative Cohen-Macaulay $R$-modules $M$
with $\cd\left(\fa,M\right)=n$.
\item[(iii)] $\text{gCM}^n_{\fa}(R)$ stands for the full subcategory of $\fa$-relative generalized Cohen-Macaulay $R$-modules
$M$ with $\cd\left(\fa,M\right)=n$.
\end{itemize}
\end{notation}

Compare the next result with \cite[Theorem 6.16]{CH}.

\begin{theorem}\label{4.1} Let $\fa$ be a proper ideal of $R$, $\RE=\underset{\fm\in \Max R}\bigoplus \RE_R(R/\fm)$ and $0\leq n\leq
c=\cd(\fa,R)$ an integer. Assume that $R$ is $\fa$-relative Cohen-Macaulay. Then for every $R$-module $M\in \text{CM}^n_{\fa}(R)$, we
have
\begin{enumerate}
\item[(i)] $\Ext_R^i(M,\Omega_{\fa})=
\begin{cases} \Hom_R(\RH_{\fa}^n(M),\RE) \hspace{.4cm}   \text{if} \hspace{.4cm}   i=c-n \\
0 \hspace{3cm}  \text{if} \hspace{.4cm}   i\neq c-n.
\end{cases}
$
\item[(ii)] $\RH_{\fa}^i\left(\Ext_R^{c-n}(M,\Omega_{\fa})\right)=
\begin{cases} \Hom_R(M,\Gamma_{\fa}(\RE)) \hspace{.4cm}  \text{if} \hspace{.4cm}   i=n \\
0 \hspace{2.9cm}   \text{if} \hspace{.4cm}   i\neq n.
\end{cases}
$
\item[(iii)] $\Ext_R^{c-n}\left(\Ext_R^{c-n}(M,\Omega_{\fa}),\Omega_{\fa}\right)\cong M\otimes_R(\underset{\fm\in \V(\fa)\cap
\Max R}\prod \widehat{R_{\fm}})$.
\end{enumerate}
\end{theorem}

\begin{prf} (i) follows by Lemma \ref{3.5}(i).

(ii) By (i), Theorem \ref{3.4}(ii) and Lemma \ref{2.6}, one has
$$\begin{array}{ll}
{\bf R}\Gamma_{\fa}(\Ext_R^{c-n}\left(M,\Omega_{\fa})\right)&\simeq {\bf R}\Gamma_{\fa}\left(\Sigma^{c-n}{\bf R}\Hom_R(M,\Omega_{\fa})\right)\\
&\simeq \Sigma^{c-n}{\bf R}\Gamma_{\fa}\left({\bf R}\Hom_R(M,\Omega_{\fa})\right)\\
&\simeq \Sigma^{c-n}{\bf R}\Hom_R\left(M,{\bf R}\Gamma_{\fa}(\Omega_{\fa})\right)\\
&\simeq  \Sigma^{c-n}{\bf R}\Hom_R\left(M,\Sigma^{-c}\Gamma_{\fa}(\RE)\right)\\
&\simeq  \Sigma^{-n}{\bf R}\Hom_R\left(M,\Gamma_{\fa}(\RE)\right)\\
&\simeq  \Sigma^{-n}\Hom_R\left(M,\Gamma_{\fa}(\RE)\right).\\
\end{array}$$
This concludes (ii).

(iii) Set $T:=\underset{\fm\in \V(\fa)\cap \Max R}\prod \widehat{R_{\fm}}$. By Theorem \ref{3.4}(i), $\id_R(\Omega_{\fa})<\infty$. Thus,
by Hom evaluation and Theorem {3.5}(iii), we deduce that $${\bf R}\Hom_R\left({\bf R}\Hom_R(M,\Omega_{\fa}),\Omega_{\fa}\right)\simeq M
\otimes_R^{\bf L}{\bf R}\Hom_R\left(\Omega_{\fa},\Omega_{\fa}\right)\simeq M\otimes_R^{\bf L}T.$$ Hence, one has
$$\begin{array}{ll}
\Ext_R^{c-n}(\Ext_R^{c-n}(M,\Omega_{\fa}),\Omega_{\fa})&\cong \RH_{-(c-n)}({\bf R}\Hom_R(\Ext_R^{c-n}(M,\Omega_{\fa}),\Omega_{\fa}))\\
&\cong \RH_{-(c-n)}({\bf R}\Hom_R\left(\Sigma^{c-n}{\bf R}\Hom_R(M,\Omega_{\fa}),\Omega_{\fa})\right)\\
&\cong \RH_{-(c-n)}(\Sigma^{-(c-n)}{\bf R}\Hom_R({\bf R}\Hom_R(M,\Omega_{\fa}),\Omega_{\fa}))\\
&\cong \RH_{0}({\bf R}\Hom_R\left({\bf R}\Hom_R(M,\Omega_{\fa}),\Omega_{\fa})\right)\\
&\cong \RH_{0}(M\otimes_R^{\bf L}T)\\
&\cong M\otimes_RT,
\end{array}$$
as desired.
\end{prf}

\begin{corollary}\label{4.2} Let $\mathfrak{J}$ denote the Jacobson radical of a complete semi-local ring $R$ and $c=\cd(\mathfrak{J},R)$.
Assume that $R$ is $\mathfrak{J}$-relative Cohen-Macaulay. Then for every integer $0\leq n\leq c$, there is a duality of categories:
\begin{displaymath}
\xymatrix{\text{CM}^n_{\mathfrak{J}}(R) \ar@<0.7ex>[rrr]^-{\Ext_R^{c-n}(-,\Omega_{\mathfrak{J}})} &
{} & {} &  \text{CM}^n_{\mathfrak{J}}(R).  \ar@<0.7ex>[lll]^-{\Ext_R^{c-n}(-,\Omega_{\mathfrak{J}})}}
\end{displaymath}
\end{corollary}

\begin{prf} Let $M$ be a nonzero finitely generated $R$-module. It is easy to see that $\cd(\mathfrak{J},M)=\dim_RM$. Although we don't
need it, it is worth mentioning that if $M$ is $\mathfrak{J}$-relative Cohen-Macaulay, then $\cd(\mathfrak{J},M)=\dim_{R_{\fm}}M_{\fm}$
for every $\fm\in \Max R\cap \Supp_RM$.

As $\dim R/\mathfrak{J}=0$, the $R$-module $\RH_{\mathfrak{J}}^c(R)$ is Artinian by \cite[Exercise 7.1.4]{BS}.
By Matlis' duality, we conclude that $\Omega_{\mathfrak{J}}$ is finitely generated, because $R$ is a complete semi-local ring. Now, the
claim follows by Theorem \ref{4.1}.
\end{prf}

One may wonder if we can replace $\text{CM}^n_{\fa}(R)$ with $\text{gCM}^n_{\fa}(R)$ in Corollary \ref{4.2}. In this regard, we
have the next result. We need the notion of finiteness dimension for its proof. Recall that for an ideal $\fa$ of $R$ and a finitely
generated $R$-module $M$, the $\fa$-{\it finiteness dimension} of $M$, denoted by $f_{\fa}(M)$, is defined as the infimum of integers
$i$ such
that the $R$-module $\RH_{\fa}^i(M)$ is not finitely generated.

\begin{theorem}\label{4.3} Let $(R,\fm)$ be a $d$-dimensional Cohen-Macaulay local ring with a dualizing module $\omega_R$. Let $M$
be an $n$-dimensional generalized Cohen-Macaulay $R$-module. Then $\Ext_R^{d-n}(M,\omega_R)$ is also a generalized Cohen-Macaulay
$R$-module of dimension $n$.
\end{theorem}

\begin{prf}  Recall that a finitely generated $R$-module $N$ is generalized Cohen-Macaulay if and only $\dim_RN\leq 0$ or $f_{\fm}(N)
=\dim_RN$. If $n\leq 0$, then $M$ is Artinian, and so $\Ext_R^{d-n}(M,\omega_R)$ is also Artinian. In this case, the claim holds
trivially. Hence, we may and do assume that $n>0$. To complete the proof, it suffices to show that $$f_{\fm}(\Ext_R^{d-n}(M,\omega_R))=n=\dim_R(\Ext_R^{d-n}(M,\omega_R)).$$

By \cite[Exercise 9.5.7]{BS}, we have $\dim R/\fq=n$ for all $\fq\in \Ass_RM\setminus \{\fm\}$ and $M_{\fp}$ is a Cohen–Macaulay
$R_{\fp}$-module for all $\fp\in \Supp_RM\setminus \{\fm\}$. Let $\fp\in \Supp_RM\setminus \{\fm\}$. As $R$ is Cohen-Macaulay and
$n>0$, it follows that $$\dim R_{\fp}-\dim_{R_{\fp}}M_{\fp}=d-n.$$ On the other hand, $(\omega_R)_{\fp}$ is the dualizing module
of the Cohen-Macaulay local ring $R_{\fp}$. Now, \cite[Theorem 3.3.10]{BH} yields that $\Ext_{R_{\fp}}^{d-n}(M_{\fp},
(\omega_R)_{\fp})$ is a Cohen-Macaulay $R_{\fp}$-module of dimension $\dim_{R_{\fp}}M_{\fp}$. In particular, $\fp\in \Supp_R(\Ext_R^{d-n}(M,\omega_R))$. Thus $$\Supp_RM=\Supp_R(\Ext_R^{d-n}(M,\omega_R)),$$ and so $\dim_R(\Ext_R^{d-n}(M,\omega_R))=n$.

For simplicity, we set $N:=\Ext_R^{d-n}(M,\omega_R)$. Then by \cite[Satz 1]{F}, we have
$$\begin{array}{ll}
f_{\fm}\left(N\right)&=\lambda_{\fm}\left(N\right)\\
&=\inf\{\depth_{\fp}N_{\fp}+\Ht((\fm+\fp)/\fp) \mid \fp\in \Supp_RN\setminus \{\fm\} \}\\
&=\inf\{\dim_{\fp}N_{\fp}+\dim R/\fp \mid \fp\in \Supp_RN\setminus \{\fm\} \}\\
&=\inf\{\dim_{\fp}M_{\fp}+\dim R/\fp \mid \fp\in \Supp_RM\setminus \{\fm\} \}\\
&=\inf\{\depth_{\fp}M_{\fp}+\dim R/\fp \mid \fp\in \Supp_RM\setminus \{\fm\} \}\\
&=\lambda_{\fm}(M)\\
&=f_{\fm}(M)\\
&=n\\
&=\dim_RN.
\end{array}$$
\end{prf}

This section's last result provides a generalizations of Grothendieck's local duality theorem.

\begin{theorem}\label{4.4} Let $\mathfrak{J}$ denote the Jacobson radical of a complete semi-local ring $R$ and $M$ a finitely generated
$R$-module. Assume that $R$ is $\mathfrak{J}$-relative Cohen-Macaulay. Then $$\RH_{\mathfrak{J}}^{i}\left(M\right)\cong \Hom_R(\Ext_R^{\dim R-i}(M,\Omega_{\mathfrak{J}}),\underset{\fm\in \Max R}\bigoplus \RE_R(R/\fm))$$ for all $i\geq 0$.
\end{theorem}

\begin{prf} Set $\RE=\underset{\fm\in \Max R}\bigoplus \RE_R(R/\fm)$ and $c=\cd(\mathfrak{J},R)$. As we saw in the proof of Corollary
\ref{4.2}, one has $\dim R=c$ and the $R$-module $\Omega_{\mathfrak{J}}$ is finitely generated.  Hence, as $\id_R(\Omega_{\mathfrak{J}})
<\infty$, it turns out that $${\bf R}\Hom_R\left(M,\Omega_{\mathfrak{J}}\right)\in \mathcal{D}^{f}_{\square}(R).$$
Theorem \ref{3.4}(iii) yields that $R\simeq {\bf R}\Hom_R(\Omega_{\mathfrak{J}},\Omega_{\mathfrak{J}})$, and so
$$M\simeq M\otimes_R^{\bf L}R\simeq M\otimes_R^{\bf L} {\bf R}\Hom_R\left(\Omega_{\mathfrak{J}},\Omega_{\mathfrak{J}}\right)\simeq
{\bf R}\Hom_R({\bf R}\Hom_R(M,\Omega_{\mathfrak{J}}),\Omega_{\mathfrak{J}}).$$ Now, Theorem \ref{3.4}(ii) implies that
${\bf R}\Gamma_{\mathfrak{J}}\left(\Omega_{\mathfrak{J}}\right)\simeq \Sigma^{-c}\RE$, and so by Lemma \ref{2.6}, we conclude that
$$\begin{array}{ll}
{\bf R}\Gamma_{\mathfrak{J}}(M)&\simeq
{\bf R}\Gamma_{\mathfrak{J}} \left({\bf R}\Hom_R({\bf R}\Hom_R(M,\Omega_{\mathfrak{J}}),\Omega_{\mathfrak{J}})\right)\\
&\simeq  {\bf R}\Hom_R\left({\bf R}\Hom_R(M,\Omega_{\mathfrak{J}}),{\bf R}\Gamma_{\mathfrak{J}}(\Omega_{\mathfrak{J}})\right)\\
&\simeq {\bf R}\Hom_R\left({\bf R}\Hom_R(M,\Omega_{\mathfrak{J}}),\Sigma^{-c}\RE\right)\\
&\simeq  \Sigma^{-c}{\bf R}\Hom_R\left({\bf R}\Hom_R(M,\Omega_{\mathfrak{J}}),\RE\right).
\end{array}$$
So for any integer $i\geq 0$, we have the following display:
$$\begin{array}{ll} \RH_{\mathfrak{J}}^{i}\left(M\right)&\cong \RH_{-i}\left({\bf R}\Gamma_{\mathfrak{J}}(M)\right)\\
&\cong \RH_{-i}\left(\Sigma^{-c}{\bf R}\Hom_R({\bf R}\Hom_R(M,\Omega_{\mathfrak{J}}),\RE)\right)\\
&\cong \RH_{c-i}\left({\bf R}\Hom_R({\bf R}\Hom_R(M,\Omega_{\mathfrak{J}}),\RE)\right)\\
&\cong \RH_{c-i}\left(\Hom_R({\bf R}\Hom_R(M,\Omega_{\mathfrak{J}}),\RE)\right)\\
&\cong \Hom_R\left(\RH_{-(c-i)}({\bf R}\Hom_R(M,\Omega_{\mathfrak{J}})),\RE\right)\\
&\cong \Hom_R\left(\Ext_R^{c-i}(M,\Omega_{\mathfrak{J}}),\RE\right).
\end{array}$$
\end{prf}

\section{Equivalences}

We need to recall the following definitions in order to present our results in this section.

\begin{definition}\label{5.1}  A finitely generated $R$-module $C$ is called {\it semidualizing} if it satisfies the following conditions:
\begin{itemize}
\item[(i)] the homothety map $\chi_C^R:R \lo \Hom_R\left(C,C\right)$ is an isomorphism, and
\item[(ii)] $\Ext^i_R\left(C,C\right)=0$ for all $i>0$.
\end{itemize}
\end{definition}

\begin{definition}\label{5.2} Let $C$ be a semidualizing module of $R$.
\begin{itemize}
\item[(i)] The Auslander class $\mathscr{A}_C\left(R\right)$ is the class of all $R$-modules $M$ for which the natural map
$\gamma_M^C:M\lo \Hom_R\left(C,C\otimes_RM\right)$ is an isomorphism, and $$\Tor^R_i\left(C,M\right)=0=\Ext_R^i
\left(C,C\otimes_RM\right)$$ for all $i\geq 1$.
\item[(ii)] The Bass class $\mathscr{B}_C\left(R\right)$ is the class of all $R$-modules $M$ for which the evaluation map
$\xi_M^C:C\otimes_R\Hom_R\left(C,M\right)\lo M$ is an isomorphism, and $$\Ext^i_R\left(C,M\right)=0=\Tor^R_i\left(C,\Hom_R\left(C,M\right)\right)$$ for all $i\geq 1$.
\end{itemize}
\end{definition}

\begin{theorem}\label{5.3}  Let $\fa$ be a proper ideal of $R$ and $C$ a semidualizing module of $R$. Then, for every non-negative
integer $n$, there is an equivalence of categories:
\begin{displaymath}
\xymatrix{\mathscr{A}_C\left(R\right)\bigcap \text{CM}^n_{\fa}(R) \ar@<0.7ex>[rrr]^-{C\otimes_R-} &
{} & {} & \mathscr{B}_C\left(R\right)\bigcap \text{CM}^n_{\fa}(R).  \ar@<0.7ex>[lll]^-{\Hom_R\left(C,-\right)}}
\end{displaymath}
\end{theorem}

\begin{prf} Because of the equivalence of categories
\begin{displaymath}
\xymatrix{\mathscr{A}_C\left(R\right) \ar@<0.7ex>[rrr]^-{C\otimes_R-} &
{} & {} & \mathscr{B}_C\left(R\right),  \ar@<0.7ex>[lll]^-{\Hom_R\left(C,-\right)}}
\end{displaymath}
for every $R$-module $L\in \mathscr{B}_C\left(R\right)$, one has $\Hom_R(C,L)\in \mathscr{A}_C\left(R\right)$, and
$L\cong C\otimes_R\Hom_R(C,L)$. Hence, to complete the proof, it is enough to show that a finitely generated $R$-module $N\in \mathscr{A}_C\left(R\right)$ belongs to $\text{CM}^n_{\fa}(R)$ if and
only if $C\otimes_RN$  belongs to $\text{CM}^n_{\fa}(R)$.

Let $M\in \mathscr{A}_C\left(R\right)$ be a finitely generated $R$-module. \cite[Lemma 3.1(ii)]{AbDT} implies that $\Supp_R\left(C\otimes_RM\right)=\Supp_RM$, and so by \cite[Theorem 2.2]{DNT}, we have $$\cd\left(\fa,M\right)=\cd\left(\fa,C\otimes_R
M\right).$$ On the other hand, \cite[Lemma 3.2(ii)]{AbDT}
implies that $$\grade\left(\fa,M\right)=\grade\left(\fa,C\otimes_RM\right).$$  Thus, $M\in \text{CM}^n_{\fa}(R)$ if and
only if $C\otimes_RM\in \text{CM}^n_{\fa}(R)$.
\end{prf}

Based on the above result, we can immediately conclude that:

\begin{corollary}\label{5.4}   Let $\fa$ be a proper ideal of $R$ and $C$ a semidualizing module of $R$. Then there is an equivalence
of categories:
\begin{displaymath}
\xymatrix{\mathscr{A}_C\left(R\right)\bigcap \text{MCM}_{\fa}(R) \ar@<0.7ex>[rrr]^-{C\otimes_R-} &
{} & {} & \mathscr{B}_C\left(R\right)\bigcap \text{MCM}_{\fa}(R).  \ar@<0.7ex>[lll]^-{\Hom_R\left(C,-\right)}}
\end{displaymath}
\end{corollary}

In the proof of the next result, we will use the notion of generalized local cohomology. Let $\fa$ be an ideal of $R$ and $M$ and
$N$ two $R$-modules. The $i$th {\it generalized local cohomology module} of $M$ and $N$ with respect to $\fa$ is defined by $$\RH^i_{\fa}(M,N):=\underset{n}{\varinjlim} \ \text{Ext}^{i}_{R}(M/\fa^{n}M,N);$$ see \cite{H}.

\begin{theorem}\label{5.5}  Let $\fa$ be a proper ideal of $R$ and $C$ a semidualizing module of $R$. Then, for every non-negative
integer $n$, there is an equivalence of categories:
\begin{displaymath}
\xymatrix{\mathscr{A}_C\left(R\right)\bigcap \text{gCM}^n_{\fa}(R) \ar@<0.7ex>[rrr]^-{C\otimes_R-} &
{} & {} & \mathscr{B}_C\left(R\right)\bigcap \text{gCM}^n_{\fa}(R).  \ar@<0.7ex>[lll]^-{\Hom_R\left(C,-\right)}}
\end{displaymath}
\end{theorem}

\begin{prf} First of all, note that a finitely generated $R$-module $M$ is $\fa$-relative generalized Cohen-Macaulay if and only if $\cd\left(\fa,M\right) \leq 0$ or $\cd\left(\fa,M\right)=\text{f}_{\fa}\left(M\right)$.

Let $M\in \mathscr{A}_C\left(R\right)$ be a finitely generated $R$-module. As we saw in the proof of Theorem \ref{5.3},
we have $$\cd\left(\fa,M\right)=\cd\left(\fa,C\otimes_RM\right).$$
	
For each pair $(L,N)$ of finitely generated $R$-modules, let $f_{\fa}(L,N)$ denote the infimum of integers $i$ such that
the $R$-module $\RH_{\fa}^i(L,N)$ is not finitely generated. Since  $\Supp_RC=\Spec R$, by
\cite[Theorem 2.1]{AsDT}, we conclude that $$f_{\fa}(C,C\otimes_RM)=f_{\fa}(R,C\otimes_RM)=f_{\fa}(C\otimes_RM).$$
As $C\otimes_RM\in \mathscr{B}_C\left(R\right)$ and $\Hom_R(C,C\otimes_RM)\cong M$, by \cite[Lemma 4.4]{AbDT}, it turns out
that $\RH_{\fa}^i(C,C\otimes_RM)\cong \RH_{\fa}^i(M)$ for all $i\geq 0$. Thus, $f_{\fa}(M)=f_{\fa}(C\otimes_RM)$.
	
Hence, for every finitely generated $R$-module $N\in \mathscr{A}_C\left(R\right)$, it follows that $N\in \text{gCM}^n_{\fa}(R)$ if and
only if $C\otimes_RN\in \text{gCM}^n_{\fa}(R)$.  Therefore, the claimed equivalence of categories follows from the following
equivalence of categories
\begin{displaymath}
\xymatrix{\mathscr{A}_C\left(R\right) \ar@<0.7ex>[rrr]^-{C\otimes_R-} &
{} & {} & \mathscr{B}_C\left(R\right).  \ar@<0.7ex>[lll]^-{\Hom_R\left(C,-\right)}}
\end{displaymath}
\end{prf}


\end{document}